\documentclass{amsart}
\usepackage{amssymb}

\theoremstyle{plain}
\newtheorem{theorem}{Theorem}
\newtheorem{axiom}[theorem]{Axiom}

\newtheorem{conjecture}[theorem]{Conjecture}
\newtheorem{corollary}[theorem]{Corollary}
\newtheorem{exercise}[theorem]{Exercise}
\newtheorem{lemma}[theorem]{Lemma}

\newtheorem{proposition}[theorem]{Proposition}

\theoremstyle{definition}
\newtheorem{definition}[theorem]{Definition}
\newtheorem{remark}[theorem]{Remark}

\newtheorem{example}[theorem]{Example}

\makeatletter
\newcount\@hour\newcount\@minute\chardef\@x10\chardef\@xv60
\def\tcitime{
\def\@time{%
  \@minute\time\@hour\@minute\divide\@hour\@xv
  \ifnum\@hour<\@x 0\fi\the\@hour:%
  \multiply\@hour\@xv\advance\@minute-\@hour
  \ifnum\@minute<\@x 0\fi\the\@minute
  }}%

\@ifundefined{hyperref}{}{}

\@ifundefined{qExtProgCall}{\def\qExtProgCall#1#2#3#4#5#6{\relax}}{}
\def\QCTOpt[#1]#2{%
  \def\QCTOptB{#1}
  \def\QCTOptA{#2}
}
\def\QCTNOpt#1{%
  \def\QCTOptA{#1}
  \let\QCTOptB\empty
}
\def\Qct{%
  \@ifnextchar[{%
    \QCTOpt}{\QCTNOpt}
}
\def\QCBOpt[#1]#2{%
  \def\QCBOptB{#1}
  \def\QCBOptA{#2}
}
\def\QCBNOpt#1{%
  \def\QCBOptA{#1}
  \let\QCBOptB\empty
}
\def\Qcb{%
  \@ifnextchar[{%
    \QCBOpt}{\QCBNOpt}
}
\def\PrepCapArgs{%
  \ifx\QCBOptA\empty
    \ifx\QCTOptA\empty
      {}%
    \else
      \ifx\QCTOptB\empty
        {\QCTOptA}%
      \else
        [\QCTOptB]{\QCTOptA}%
      \fi
    \fi
  \else
    \ifx\QCBOptA\empty
      {}%
    \else
      \ifx\QCBOptB\empty
        {\QCBOptA}%
      \else
        [\QCBOptB]{\QCBOptA}%
      \fi
    \fi
  \fi
}
\newcount\GRAPHICSTYPE
\GRAPHICSTYPE=\z@
\def\GRAPHICSPS#1{%
 \ifcase\GRAPHICSTYPE
   \special{ps: #1}%
 \or
   \special{language "PS", include "#1"}%
 \fi
}%
\def\graffile#1#2#3#4{%
    \leavevmode
    \raise -#4 \BOXTHEFRAME{%
        \hbox to #2{\raise #3\hbox to #2{\null #1\hfil}}}%
}%
\def\draftbox#1#2#3#4{%
 \leavevmode\raise -#4 \hbox{%
  \frame{\rlap{\protect\tiny #1}\hbox to #2%
   {\vrule height#3 width\z@ depth\z@\hfil}%
  }%
 }%
}%
\newcount\draft
\draft=\z@

\newif\ifwasdraft
\wasdraftfalse

\def\GRAPHIC#1#2#3#4#5{%
 \ifnum\draft=\@ne\draftbox{#2}{#3}{#4}{#5}%
  \else\graffile{#1}{#3}{#4}{#5}%
  \fi
 }%
\def\addtoLaTeXparams#1{%
    \edef\LaTeXparams{\LaTeXparams #1}}%
\newif\ifBoxFrame \BoxFramefalse
\newif\ifOverFrame \OverFramefalse
\newif\ifUnderFrame \UnderFramefalse

\def\BOXTHEFRAME#1{%
   \hbox{%
      \ifBoxFrame
         \frame{#1}%
      \else
         {#1}%
      \fi
   }%
}

\def\doFRAMEparams#1{\BoxFramefalse\OverFramefalse\UnderFramefalse\readFRAMEparams#1\end}%
\def\readFRAMEparams#1{%
 \ifx#1\end%
  \let\next=\relax
  \else
  \ifx#1i\dispkind=\z@\fi
  \ifx#1d\dispkind=\@ne\fi
  \ifx#1f\dispkind=\tw@\fi
  \ifx#1t\addtoLaTeXparams{t}\fi
  \ifx#1b\addtoLaTeXparams{b}\fi
  \ifx#1p\addtoLaTeXparams{p}\fi
  \ifx#1h\addtoLaTeXparams{h}\fi
  \ifx#1X\BoxFrametrue\fi
  \ifx#1O\OverFrametrue\fi
  \ifx#1U\UnderFrametrue\fi
  \ifx#1w
    \ifnum\draft=1\wasdrafttrue\else\wasdraftfalse\fi
    \draft=\@ne
  \fi
  \let\next=\readFRAMEparams
  \fi
 \next
 }%
\def\IFRAME#1#2#3#4#5#6{%
      \bgroup
      \let\QCTOptA\empty
      \let\QCTOptB\empty
      \let\QCBOptA\empty
      \let\QCBOptB\empty
      #6%
      \parindent=0pt%
      \leftskip=0pt
      \rightskip=0pt
      \setbox0 = \hbox{\QCBOptA}%
      \@tempdima = #1\relax
      \ifOverFrame
          \typeout{This is not implemented yet}%
          \show\HELP
      \else
         \ifdim\wd0>\@tempdima
            \advance\@tempdima by \@tempdima
            \ifdim\wd0 >\@tempdima
               \textwidth=\@tempdima
               \setbox1 =\vbox{%
                  \noindent\hbox to \@tempdima{\hfill\GRAPHIC{#5}{#4}{#1}{#2}{#3}\hfill}\\%
                  \noindent\hbox to \@tempdima{\parbox[b]{\@tempdima}{\QCBOptA}}%
               }%
               \wd1=\@tempdima
            \else
               \textwidth=\wd0
               \setbox1 =\vbox{%
                 \noindent\hbox to \wd0{\hfill\GRAPHIC{#5}{#4}{#1}{#2}{#3}\hfill}\\%
                 \noindent\hbox{\QCBOptA}%
               }%
               \wd1=\wd0
            \fi
         \else
            \ifdim\wd0>0pt
              \hsize=\@tempdima
              \setbox1 =\vbox{%
                \unskip\GRAPHIC{#5}{#4}{#1}{#2}{0pt}%
                \break
                \unskip\hbox to \@tempdima{\hfill \QCBOptA\hfill}%
              }%
              \wd1=\@tempdima
           \else
              \hsize=\@tempdima
              \setbox1 =\vbox{%
                \unskip\GRAPHIC{#5}{#4}{#1}{#2}{0pt}%
              }%
              \wd1=\@tempdima
           \fi
         \fi
         \@tempdimb=\ht1
         \advance\@tempdimb by \dp1
         \advance\@tempdimb by -#2%
         \advance\@tempdimb by #3%
         \leavevmode
         \raise -\@tempdimb \hbox{\box1}%
      \fi
      \egroup%
}%
\def\DFRAME#1#2#3#4#5{%
 \begin{center}
     \let\QCTOptA\empty
     \let\QCTOptB\empty
     \let\QCBOptA\empty
     \let\QCBOptB\empty
     \ifOverFrame 
        #5\QCTOptA\par
     \fi
     \GRAPHIC{#4}{#3}{#1}{#2}{\z@}
     \ifUnderFrame 
        \nobreak\par #5\QCBOptA
     \fi
 \end{center}%
 }%
\def\FFRAME#1#2#3#4#5#6#7{%
 \begin{figure}[#1]%
  \let\QCTOptA\empty
  \let\QCTOptB\empty
  \let\QCBOptA\empty
  \let\QCBOptB\empty
  \ifOverFrame
    #4
    \ifx\QCTOptA\empty
    \else
      \ifx\QCTOptB\empty
        \caption{\QCTOptA}%
      \else
        \caption[\QCTOptB]{\QCTOptA}%
      \fi
    \fi
    \ifUnderFrame\else
      \label{#5}%
    \fi
  \else
    \UnderFrametrue%
  \fi
  \begin{center}\GRAPHIC{#7}{#6}{#2}{#3}{\z@}\end{center}%
  \ifUnderFrame
    #4
    \ifx\QCBOptA\empty
      \caption{}%
    \else
      \ifx\QCBOptB\empty
        \caption{\QCBOptA}%
      \else
        \caption[\QCBOptB]{\QCBOptA}%
      \fi
    \fi
    \label{#5}%
  \fi
  \end{figure}%
 }%
\newcount\dispkind%

\def\makeactives{
  \catcode`\"=\active
  \catcode`\;=\active
  \catcode`\:=\active
  \catcode`\'=\active
  \catcode`\~=\active
}
\bgroup
   \makeactives
   \gdef\activesoff{%
      \def"{\string"}
      \def;{\string;}
      \def:{\string:}
      \def'{\string'}
      \def~{\string~}
    }
\egroup

\def\FRAME#1#2#3#4#5#6#7#8{%
 \bgroup
 \@ifundefined{bbl@deactivate}{}{\activesoff}
 \ifnum\draft=\@ne
   \wasdrafttrue
 \else
   \wasdraftfalse%
 \fi
 \def\LaTeXparams{}%
 \dispkind=\z@
 \def\LaTeXparams{}%
 \doFRAMEparams{#1}%
 \ifnum\dispkind=\z@\IFRAME{#2}{#3}{#4}{#7}{#8}{#5}\else
  \ifnum\dispkind=\@ne\DFRAME{#2}{#3}{#7}{#8}{#5}\else
   \ifnum\dispkind=\tw@
    \edef\@tempa{\noexpand\FFRAME{\LaTeXparams}}%
    \@tempa{#2}{#3}{#5}{#6}{#7}{#8}%
    \fi
   \fi
  \fi
  \ifwasdraft\draft=1\else\draft=0\fi{}%
  \egroup
 }%
\def\TEXUX#1{"texux"}

\def\limfunc#1{\mathop{\rm #1}}%

\long\def\QQQ#1#2{%
     \long\expandafter\def\csname#1\endcsname{#2}}%
\@ifundefined{QTP}{\def\QTP#1{}}{}
\@ifundefined{QEXCLUDE}{\def\QEXCLUDE#1{}}{}
\@ifundefined{Qlb}{}{}
\@ifundefined{Qlt}{}{}
\long\def\QQA#1#2{}%
\def\QTR#1#2{{\csname#1\endcsname #2}}
\def\EXPAND#1[#2]#3{}%
\def\NOEXPAND#1[#2]#3{}%
\def\LaTeXparent#1{}%
\def\ChildStyles#1{}%
\def\ChildDefaults#1{}%
\def\QTagDef#1#2#3{}%
\@ifundefined{StyleEditBeginDoc}{}{}
\def\QQfnmark#1{\footnotemark}

\@ifundefined{INDEX}{\def\INDEX#1#2{}{}}{}%
\@ifundefined{SUBINDEX}{\def\SUBINDEX#1#2#3{}{}{}}{}%
\@ifundefined{initial}%
   {\def\initial#1{\bigbreak{\raggedright\large\bf #1}\kern 2\p@\penalty3000}}%
   {}%
\@ifundefined{entry}{}{}%
\@ifundefined{primary}{}{}%
\@ifundefined{secondary}{}{}%
\@ifundefined{ZZZ}{}{\makeatletter\input gnuindex.sty\makeatother\makeindex\makeatletter}%
\@ifundefined{abstract}{%
 \def\abstract{%
  \if@twocolumn
   \section*{Abstract (Not appropriate in this style!)}%
   \else \small 
   \begin{center}{\bf Abstract\vspace{-.5em}\vspace{\z@}}\end{center}%
   \quotation 
   \fi
  }%
 }{%
 }%
\@ifundefined{endabstract}{\def\endabstract
  {\if@twocolumn\else\endquotation\fi}}{}%
\@ifundefined{maketitle}{\def\maketitle#1{}}{}%
\@ifundefined{affiliation}{\def\affiliation#1{}}{}%
\@ifundefined{proof}{}{}%
\@ifundefined{endproof}{}{}%
\@ifundefined{newfield}{\def\newfield#1#2{}}{}%
\@ifundefined{chapter}{\def\chapter#1{\par(Chapter head:)#1\par }%
 \newcount\c@chapter}{}%
\@ifundefined{part}{\def\part#1{\par(Part head:)#1\par }}{}%
\@ifundefined{section}{\def\section#1{\par(Section head:)#1\par }}{}%
\@ifundefined{subsection}{\def\subsection#1%
 {\par(Subsection head:)#1\par }}{}%
\@ifundefined{subsubsection}{\def\subsubsection#1%
 {\par(Subsubsection head:)#1\par }}{}%
\@ifundefined{paragraph}{\def\paragraph#1%
 {\par(Subsubsubsection head:)#1\par }}{}%
\@ifundefined{subparagraph}{\def\subparagraph#1%
 {\par(Subsubsubsubsection head:)#1\par }}{}%
\@ifundefined{therefore}{}{}%
\@ifundefined{backepsilon}{}{}%
\@ifundefined{yen}{}{}%
\@ifundefined{registered}{%
   \def\registered{\relax\ifmmode{}\r@gistered
                    \else$\m@th\r@gistered$\fi}%
 \def\r@gistered{^{\ooalign
  {\hfil\raise.07ex\hbox{$\scriptstyle\rm\text{R}$}\hfil\crcr
  \mathhexbox20D}}}}{}%
\@ifundefined{Eth}{}{}%
\@ifundefined{eth}{}{}%
\@ifundefined{Thorn}{}{}%
\@ifundefined{thorn}{}{}%
\@ifundefined{degree}{}{}%
\newdimen\theight
\def\Column{%
 \vadjust{\setbox\z@=\hbox{\scriptsize\quad\quad tcol}%
  \theight=\ht\z@\advance\theight by \dp\z@\advance\theight by \lineskip
  \kern -\theight \vbox to \theight{%
   \rightline{\rlap{\box\z@}}%
   \vss
   }%
  }%
 }%
\def\qed{%
 \ifhmode\unskip\nobreak\fi\ifmmode\ifinner\else\hskip5\p@\fi\fi
 \hbox{\hskip5\p@\vrule width4\p@ height6\p@ depth1.5\p@\hskip\p@}%
 }%
\def\miss{\hbox{\vrule height2\p@ width 2\p@ depth\z@}}%
\def\tcol#1{{\baselineskip=6\p@ \vcenter{#1}} \Column}  %
\def\newfmtname{LaTeX2e}
\def\chkcompat{%
   \if@compatibility
   \else
     \usepackage{latexsym}
   \fi
}

\ifx\fmtname\newfmtname
  \DeclareOldFontCommand{\rm}{\normalfont\rmfamily}{\mathrm}
  \DeclareOldFontCommand{\sf}{\normalfont\sffamily}{\mathsf}
  \DeclareOldFontCommand{\tt}{\normalfont\ttfamily}{\mathtt}
  \DeclareOldFontCommand{\bf}{\normalfont\bfseries}{\mathbf}
  \DeclareOldFontCommand{\it}{\normalfont\itshape}{\mathit}
  \DeclareOldFontCommand{\sl}{\normalfont\slshape}{\@nomath\sl}
  \DeclareOldFontCommand{\sc}{\normalfont\scshape}{\@nomath\sc}
  \chkcompat
\fi

\def\alpha{{\Greekmath 010B}}%
\def\beta{{\Greekmath 010C}}%
\def\gamma{{\Greekmath 010D}}%
\def\delta{{\Greekmath 010E}}%
\def\epsilon{{\Greekmath 010F}}%
\def\zeta{{\Greekmath 0110}}%
\def\eta{{\Greekmath 0111}}%
\def\theta{{\Greekmath 0112}}%
\def\iota{{\Greekmath 0113}}%
\def\kappa{{\Greekmath 0114}}%
\def\lambda{{\Greekmath 0115}}%
\def\mu{{\Greekmath 0116}}%
\def\nu{{\Greekmath 0117}}%
\def\xi{{\Greekmath 0118}}%
\def\pi{{\Greekmath 0119}}%
\def\rho{{\Greekmath 011A}}%
\def\sigma{{\Greekmath 011B}}%
\def\tau{{\Greekmath 011C}}%
\def\upsilon{{\Greekmath 011D}}%
\def\phi{{\Greekmath 011E}}%
\def\chi{{\Greekmath 011F}}%
\def\psi{{\Greekmath 0120}}%
\def\omega{{\Greekmath 0121}}%
\def\varepsilon{{\Greekmath 0122}}%
\def\vartheta{{\Greekmath 0123}}%
\def\varpi{{\Greekmath 0124}}%
\def\varrho{{\Greekmath 0125}}%
\def\varsigma{{\Greekmath 0126}}%
\def\varphi{{\Greekmath 0127}}%

\def\nabla{{\Greekmath 0272}}
\def\FindBoldGroup{%
   {\setbox0=\hbox{$\mathbf{x\global\edef\theboldgroup{\the\mathgroup}}$}}%
}

\def\Greekmath#1#2#3#4{%
    \if@compatibility
        \ifnum\mathgroup=\symbold
           \mathchoice{\mbox{\boldmath$\displaystyle\mathchar"#1#2#3#4$}}%
                      {\mbox{\boldmath$\textstyle\mathchar"#1#2#3#4$}}%
                      {\mbox{\boldmath$\scriptstyle\mathchar"#1#2#3#4$}}%
                      {\mbox{\boldmath$\scriptscriptstyle\mathchar"#1#2#3#4$}}%
        \else
           \mathchar"#1#2#3#4%
        \fi 
    \else 
        \FindBoldGroup
        \ifnum\mathgroup=\theboldgroup 
           \mathchoice{\mbox{\boldmath$\displaystyle\mathchar"#1#2#3#4$}}%
                      {\mbox{\boldmath$\textstyle\mathchar"#1#2#3#4$}}%
                      {\mbox{\boldmath$\scriptstyle\mathchar"#1#2#3#4$}}%
                      {\mbox{\boldmath$\scriptscriptstyle\mathchar"#1#2#3#4$}}%
        \else
           \mathchar"#1#2#3#4%
        \fi     	    
	  \fi}

\newif\ifGreekBold  \GreekBoldfalse
\let\SAVEPBF=\pbf
\def\pbf{\GreekBoldtrue\SAVEPBF}%

\@ifundefined{theorem}{\newtheorem{theorem}{Theorem}}{}
\@ifundefined{lemma}{\newtheorem{lemma}[theorem]{Lemma}}{}
\@ifundefined{corollary}{}{}
\@ifundefined{conjecture}{}{}
\@ifundefined{proposition}{}{}
\@ifundefined{axiom}{}{}
\@ifundefined{remark}{}{}
\@ifundefined{example}{}{}
\@ifundefined{exercise}{}{}
\@ifundefined{definition}{\newtheorem{definition}{Definition}}{}

\@ifundefined{mathletters}{%
  \newcounter{equationnumber}  
  \def\mathletters{%
     \addtocounter{equation}{1}
     \edef\@currentlabel{\theequation}%
     \setcounter{equationnumber}{\c@equation}
     \setcounter{equation}{0}%
     \edef\theequation{\@currentlabel\noexpand\alph{equation}}%
  }
  
}{}

\@ifundefined{BibTeX}{%
    \def\BibTeX{{\rm B\kern-.05em{\sc i\kern-.025em b}\kern-.08em
                 T\kern-.1667em\lower.7ex\hbox{E}\kern-.125emX}}}{}%
\@ifundefined{AmS}%
    {\def\AmS{{\protect\usefont{OMS}{cmsy}{m}{n}%
                A\kern-.1667em\lower.5ex\hbox{M}\kern-.125emS}}}{}%
\@ifundefined{AmSTeX}{}{}%
\begin{document}
\title{On Whitehead precovers}
\author{Paul C. Eklof}
\thanks{First author partially supported by NSF DMS 98-03126.}
\address[Eklof]{Math Dept, UCI\\
Irvine, CA 92697-3875}
\author{Saharon Shelah}
\thanks{Second author supported  by
the German-Israeli Foundation for Scientific Research \& Development. Publication
749.}
\address[Shelah]{Institute of Mathematics, Hebrew University\\
Jerusalem 91904, Israel}
\date{\today}
\maketitle

\begin{abstract}
It is proved undecidable in ZFC + GCH whether every $\Bbb{Z}$-module has a $%
^{\perp }\{\Bbb{Z\}}$-precover.
\end{abstract}


Let $\mathcal{F}$ be a class of $R$-modules of the form 
\[
^{\perp }\mathcal{C}=\{A:\hbox{Ext}(A,C)=0\hbox{ for all }C\in \mathcal{C}\}
\]
for some class $\mathcal{C}$. The first author and Jan Trlifaj proved \cite
{ET1} that a sufficient condition for every module $M$ to have an $\mathcal{F%
}$-precover is that there is a module $B$ such that $\mathcal{F}^{\perp
}=\{B\}^{\perp }$ ($=\{A:\hbox{Ext}(B,A)=0\}$). In \cite{ET2}, generalizing
a method used by Enochs \cite{BBE} to prove the Flat Cover Conjecture, it is
proved that this sufficient condition holds whenever $\mathcal{C}$ is a
class of pure-injective modules; moreover, for $R$ a Dedekind domain, the
sufficient condition holds whenever $\mathcal{C}$ is a class of cotorsion
modules. The following is also proved in \cite{ET2}:

\begin{theorem}
\label{hered}It is consistent with ZFC + GCH that for any hereditary ring $R$
and any $R$-module $N$, there is an $R$-module $B$ such that $(^{\perp
}\{N\})^{\perp }=\{B\}^{\perp }$ and hence every $R$-module has a $^{\perp
}\{N\}$-precover.
\end{theorem}

This is a generalization of a result proved by the second author for the
class $\mathcal{W}$ of all Whitehead groups ( $=$ $^{\perp }\{\Bbb{Z\}}$):

\begin{theorem}
\label{free}It is consistent with ZFC + GCH that $\mathcal{W}^{\perp
}=\{B\}^{\perp }$ where $B$ is any free abelian group.
\end{theorem}

\noindent \textsc{Proof}. The second author proved that G\"{o}del's Axiom of
Constructibility (V = L) implies that $\mathcal{W}$ is exactly the class of
free groups. (See \cite{Sh74} or \cite{E76}.) Under this hypothesis (which
implies GCH), $\mathcal{W}^{\perp }$ is the class of all groups; if we take $%
B$ to be any free group, then $\{B\}^{\perp }$ is also the class of all
groups.\qed  

\medskip 

Our main result here is that the conclusions of Theorem \ref{hered} are not
provable in ZFC + GCH for $N=\Bbb{Z}=R$:

\begin{theorem}
\label{noncov}It is consistent with ZFC + GCH that $\Bbb{Q}$ does not have a 
$\mathcal{W}$-precover.
\end{theorem}

An immediate consequence is:

\begin{theorem}
\label{noB}It is consistent with ZFC + GCH that there is no abelian group $B$
such that $\mathcal{W}^{\perp }=\{B\}^{\perp }$.
\end{theorem}

Theorem \ref{noncov} follows easily from the following:

\begin{theorem}
\label{noncov2}It is consistent with ZFC + GCH that for every Whitehead
group $B$ there is an uncountable Whitehead group $G=G_{B}$ such that every
homomorphism from $G$ to $B$ has finitely-generated range.
\end{theorem}

\noindent \textbf{Proof of Theorem \ref{noncov} from Theorem \ref{noncov2}.}
Suppose that $f:B\rightarrow \Bbb{Q}$ is a $\mathcal{W}$-precover of $\Bbb{Q}
$. Let $G$ be as in Theorem \ref{noncov2} for this $B$. Since $\Bbb{Q}$ is
injective and $G$ has infinite rank, there is a surjective homomorphism $%
g:G\rightarrow \Bbb{Q}$. But then clearly there is no $h:G\rightarrow B$
such that $f\circ h=g.$ \qed  

\medskip 

We get the hypothesis of Theorem \ref{noncov2} from the following:

\begin{theorem}
\label{noncov3}Assume GCH. Suppose that for every Whitehead group $A$ of
infinite rank, there is a Whitehead group $H_{A}$ of cardinality $\leq
|A|^{+}$ such that $\limfunc{Ext}(H_{A},A)\neq 0$. Then for every Whitehead
group $B$ there is an uncountable Whitehead group $G$ such that every
homomorphism from $G$ to $B$ has finitely-generated range.
\end{theorem}

\noindent \textsc{Proof}. Let $\lambda =\mu ^{+}$ where $\mu >|B|+\aleph
_{1} $. Then $\diamondsuit _{\lambda }$ holds, by GCH, and we will use it to
construct the group structure on a set $G$ of size $\lambda $. We can write $%
G=\bigcup_{\nu <\lambda }G_{\nu }$ as the union of a continuous chain of
sets such that for all $\nu <\lambda $, $|G_{\nu +1}-G_{\nu }|=\mu $. Now $%
\diamondsuit _{\lambda }$ gives us a family $\{h_{\nu }:\nu \in \lambda \}$
of functions $h_{\nu }:G_{\nu }\rightarrow B$ such that for every function $%
f:G\rightarrow B$, $\{\nu \in \lambda :f\upharpoonright G_{\nu }=h_{\nu }\}$
is stationary.

Suppose that the group structure on $G_{\nu }$ has been defined and consider 
$h_{\nu }$; if the range of $h_{\nu }$ is of finite rank, define the group
structure on $G_{\nu +1}$ in any way which extends that on $G_{\nu }$.
Otherwise, let $A$ be the range of $h_{\nu }$ and let $H_{A}$ be as in the
hypothesis. Without loss of generality, $|H_{A}|=\mu $; write $H_{A}=F/K$
where $F$ is a free group of rank $\mu $. By a standard homological
argument, there is a homomorphism $\psi :K\rightarrow A$ which does not
extend to a homomorphism $\varphi :F\rightarrow A$. Since $K$ is free and $%
h_{\nu }:G_{\nu }\rightarrow B$ is onto $A$, there is a homomorphism $\theta
:K\rightarrow G_{\nu }$ such that $h_{\nu }\circ \theta =\psi $. Now form
the pushout 
\[
\begin{array}{lll}
F & \rightarrow  & G_{\nu +1} \\ 
\uparrow  &  & \uparrow  \\ 
K & \stackrel{\theta }{\rightarrow } & G_{\nu }
\end{array}
\]
to define the group structure on $G_{\nu +1}$ (cf. \cite[proof of Theorem 2]
{ET1}). Then $G_{\nu +1}/G_{\nu }\cong F/K\cong H_{A}$ so it is Whitehead.
Moreover, $h_{\nu }$ does not extend to a homomorphism from $F$  into $A$,
else $\psi $ does. This completes the definition of $G$. Notice that $G$ is
a Whitehead group since all quotients $G_{\nu +1}/G_{\nu }$ are Whitehead
(cf. \cite[Lemma 1]{ET1}).

Now given any homomorphism $f:G\rightarrow B$, let $A\subseteq B$ be the
range of $f$. Since $|A|<|G|=\lambda $, $\{\nu \in \lambda :f[G_{\nu }]=A\}$
is a club in $\lambda $; hence there exists $\nu \in \lambda $ such that $%
f\upharpoonright G_{\nu }=h_{\nu }$ and the range of $h_{\nu }$ is $A$. If $A
$ is of infinite rank, we have constructed $G_{\nu +1}$ so that $%
f\upharpoonright G_{\nu }$ does not extend to $G$, which is a contradiction.
So we must conclude that the range of $f$ is of finite rank. \qed  

\medskip

Now our main task is to show that there is a model of ZFC + GCH where the
hypothesis of Theorem \ref{noncov3} hold. As a warm-up exercise, however, we
will begin in the next section with a direct proof of Theorem \ref{noB};
this is equivalent to the consistency of a weaker assumption than the
hypothesis of Theorem \ref{noncov3}.

\section{$\mathcal{W}$ is not cogenerated by a set}

Theorem \ref{noB} is equivalent to the statement that it is consistent with
ZFC + GCH that for every W-group $B$ we can find a W-group $A\in
\{B\}^{\perp }$ such that there is a W-group $H_{A}$ with $\limfunc{Ext}%
(H_{A},A)\neq 0$. The proof will use the following consequence of Theorem 2
of \cite{ET1}:

\begin{theorem}
\label{vanext}Let $\mu $ be a  cardinal  $>\kappa $ such that $\mu ^{\kappa
}=\kappa $ and let $B$ be a group of cardinality $\leq \kappa $. Then there
is a group $A\in \{B\}^{\perp }$ such that $A=\cup _{\nu <\mu }A_{\nu
}^{\prime }$ (continuous), $A_{0}^{\prime }=0$, and such that for all $\nu
<\sigma $, $A/A_{\nu }^{\prime }$ is isomorphic to $B$. \qed  
\end{theorem}

\noindent 

\noindent \textbf{\ Proof of Theorem \ref{noB}. }We will use the fact that
the following principle is consistent with ZFC + GCH (cf. \cite{ES}):

\begin{quote}
\textbf{(UP)} For every cardinal $\sigma $ of the form $\tau ^{+}$ where $%
\tau $ is singular of cofinality $\omega $ there is a stationary subset $S$
of $\sigma $ consisting of limit ordinals of cofinality $\omega $ and a
ladder system $\bar{\zeta}=\{\zeta _{\delta }:\delta \in S\}$ which has the $%
\omega $-uniformization property, that is, for every family $\{c_{\delta
}:\delta \in S\}$ of functions $c_{\delta }:\omega \rightarrow \omega $,
there is a function $h:\sigma \rightarrow \omega $ such that for every $%
\delta \in S$, $h(\zeta _{\delta }(n))=c_{\delta }(n)$ for almost all $n\in
\omega $.
\end{quote}

\noindent We work in a model of GCH plus UP. \textbf{\ }Let $\kappa =|B|$
and let $\mu $ be a singular cardinal of cofinality $\sigma >\kappa $ such
that $\sigma $ is the successor of a singular cardinal of cofinality $\omega 
$. Then $\mu ^{\kappa }=\kappa $.  Let $A=\cup _{\nu <\sigma }A_{\nu
}^{\prime }$ be as in Theorem \ref{vanext} for this $B$ and $\mu $. Choose a
strictly increasing continuous function $\xi :\sigma \rightarrow \mu $ whose
range in cofinal in $\mu $ and let $A_{\nu }=A_{\xi (\nu )}^{\prime }$. Let $%
\bar{\zeta}=\{\zeta _{\delta }:\delta \in S\}$ be as in (UP).

Let $H_{A}=F/K$ where $F$ is the free group on symbols $\{y_{\delta ,n}:$ $%
\delta \in S$, $n\in \omega \}\cup \{x_{j}:$ $j<\sigma \}$ and $K$ is the
subgroup with basis $\{w_{\delta ,n}:$ $\delta \in S$, $n\in \omega \}$
where 
\begin{equation}
w_{\delta ,n}=2y_{\delta ,n+1}-y_{\delta ,n}+x_{\zeta _{\delta }(n)}.
\label{eqnsG}
\end{equation}

Then $H_{A}$ is a group of cardinality $\sigma $ and the $\omega $%
-uniformization property of $\bar{\zeta}$ implies that $H_{A}$ is a
Whitehead group (see \cite[XII.3]{EM} or  \cite{T}).

Now for all $\nu <\mu $, $A/A_{\nu }$ is a W-group and hence strongly $%
\aleph _{1}$-free, since CH holds, so it has a pure free subgroup $C/A_{\nu }
$ of rank $\omega $ with basis $\{t_{\nu ,n}+A_{\delta }:n\in \omega \}$
such that $A/C$ is $\aleph _{1}$-free. Then $a_{\delta }=\Sigma _{n\in
\omega }2^{n}(t_{\nu ,n}+A_{\delta })$ is in the $2$-adic completion of $%
A/A_{\delta }$ but not in $A/A_{\delta }$.

Now define $\psi :K\rightarrow A$ such that $\psi (w_{\delta ,n})=t_{\delta
,n}$ for all $\delta \in S$, $n\in \omega $. We claim that $\psi $ does not
extend to a homomorphism $\varphi :F\rightarrow A$. Suppose, to the
contrary, that it does. The set of $\delta <\sigma $ such that $\varphi
(x_{j})\in A_{\delta }$ for all $j<\delta $ is a club, $C$, in $\sigma $, so
there exists $\delta \in S\cap C$. We will contradict the choice of $%
a_{\delta }$ for this $\delta $.

We work in $A/A_{\delta }$. Let $c_{n}=\varphi (y_{\delta ,n})+A_{\delta }$.
Then by applying $\varphi $ to the equations (\ref{eqnsG}) and since $%
\varphi (x_{j})\in A_{\delta }$ for all $j<\delta $ we have that for all $%
n\in \omega $, 
\[
t_{\delta ,n}+A_{\delta }=2c_{n+1}+c_{n}.
\]
It follows that $a_{\delta }=c_{0}$ is in $A/A_{\delta }$, a contradiction. 
\qed  

\medskip

This completes the proof of the weaker Theorem \ref{noB}. In Theorem \ref
{noncov3}, we are not able to choose the Whitehead group $A$, but must find
an $H_{A}$ for every $A$. In the next section we discuss ways to insure that 
$\limfunc{Ext}(H_{A},A)$ is non-zero, and in the following section we deal
with how to make $H_{A}$ a Whitehead group, and then finish the proof of the
main theorem.

.

\section{How to make Ext not vanish}

We begin by proving some general properties of decompositions of Whitehead
groups assuming GCH. We use the result of Gregory and Shelah (cf. \cite{Greg}%
, \cite{ShL(Q)}) that GCH implies $\diamondsuit _{\lambda }$ for every
successor cardinal $\lambda >\aleph _{1}$, and the result of Devlin and
Shelah \cite{DS} that CH implies weak diamond, $\Phi _{\aleph _{1}}$,  at $%
\aleph _{1}$. We will also make repeated use of the fact (cf. \cite{Sh74}, 
\cite[Chap XII]{EM}) that if $A=\bigcup_{\alpha <\lambda }A_{\alpha }$ is a $%
\lambda $ filtration of a group of cardinality $\lambda $ and if $%
\diamondsuit _{\lambda }(E)$ holds where $E=\{\alpha \in \lambda :\exists
\beta >\alpha $ s.t. $A_{\beta }/A_{\alpha }$ is not Whitehead$\}$, then $A$
is not a Whitehead group.

\begin{lemma}
\label{nice}Let $A$ be a Whitehead group of cardinality $\lambda =\mu ^{+}$
and write $A=\bigcup_{\alpha <\lambda }A_{\alpha }$ as the continuous union
of a chain of subgroups of cardinality $\mu $. Let $S(A)\stackrel{\text{def}%
}{=}$ $\{\alpha \in \lambda :\alpha $ is a limit and $A_{\tau }/A_{\alpha }$
is Whitehead for all $\tau >\alpha \}$. If $\diamondsuit _{\lambda }(Y)$
holds for some subset $Y$ of $\lambda $, then $Y\cap S(A)$ is stationary. In
particular, assuming GCH, $S(A)$ is stationary.
\end{lemma}

\noindent \textsc{Proof}. Suppose $Y\cap S(A)$ is not stationary in $\lambda 
$, and let $C$ be a club in its complement; then we can define a continuous
increasing function $f:\lambda \rightarrow C$ such that for all $\alpha \in
\lambda $, if $f(\alpha )\in Y$, then $A_{f(\alpha +1)}/A_{f(\alpha )}$ is
not Whitehead. But then $\diamondsuit _{\lambda }(Y\cap \limfunc{im}(f))$
holds and implies that $A=\bigcup_{\alpha \in \lambda }A_{f(\alpha )}$ is
not Whitehead. \qed  

\medskip 

We can say, for short, that $A/A_{\alpha }$ is \textit{locally Whitehead }%
when   $\alpha \in S(A)$,  since every subgroup of $A/A_{\alpha }$ of
cardinality $<\lambda $ is Whitehead.

\begin{lemma}
\label{nice2}Assume GCH. Let $A$ be a Whitehead group of cardinality $\mu $
(possibly a singular cardinal). Then we can write $A=\bigcup_{\nu <\mu
}A_{\nu }$ as the continuous union of a chain of subgroups of cardinality $%
<\mu $ such that for all $\nu <\mu $, $A/A_{\nu +1}$ is $\aleph _{1}$-free.
\end{lemma}

\noindent \textsc{Proof}. If suffices to show that every subgroup $X$ of $A$
of cardinality $\kappa $ $<\mu $ is contained in a subgroup $N$ of
cardinality $\kappa $ such that $N^{\prime }/N$ is free whenever $N\subseteq
N^{\prime }\subseteq A$ and $N^{\prime }/N$ is countable. But if $X$ is a
counterexample, then we can build a chain $\{N_{\alpha }:\alpha <\kappa
^{+}\}$ such that $N_{0}=X$ and for all $\alpha <\kappa ^{+}$, $N_{\alpha
+1}/N_{\alpha }$ is countable and not free, and hence is not Whitehead. We
obtain a contradiction since then $\diamondsuit _{\kappa ^{+}}$ implies that 
$\bigcup_{\alpha <\kappa ^{+}}N_{\alpha }$ is not Whitehead. \qed  

\medskip

We now give sufficient conditions for $\limfunc{Ext}(H,A)$ to be non-zero,
when $H$ is given by a relative simple set of relations defined using ladder
systems (see the definition below). The analysis will be divided into cases,
depending on whether the cardinality of $A$ is singular, the successor of a
regular cardinal, or the successor of a singular cardinal.

The following concrete description of a group is in the spirit of the
general constructions in, for example, \cite{T} or \cite[XII.3.4]{EM} but is
a little more complicated since it is ``two step'': involving a system of
ladders of length $\limfunc{cf}(\mu )$ and another system of ladders of
length $\omega $ (if $\limfunc{cf}(\mu )>\aleph _{0}$).

\begin{definition}
\medskip \label{built}Let $\mu $ be a cardinal of cofinality $\sigma $ ($%
\leq \mu $). Let $S$ be a subset of $\lambda =\mu ^{+}$ consisting of
ordinals of cofinality $\sigma $ and $\bar{\eta}=\{\eta _{\delta }:\delta
\in S\}$ a ladder system on $S$, that is, a family of functions $\eta
_{\delta }:\sigma \rightarrow \sigma $ which are strictly increasing and
cofinal. If $\sigma >\aleph _{0}$, let $E$ be a stationary subset of $\sigma 
$ consisting of limit ordinals of cofinality $\omega $ and let $\bar{\zeta}%
=\{\zeta _{\nu }:\nu \in E\}$ be a ladder system on $E$. We will say that $H$
is \textit{the group built on }$\bar{\eta}$\textit{\ and }$\bar{\zeta}$ if $%
H\cong F/K$ where $F$ is the free group on symbols $\{y_{\delta ,\nu
,n}:\delta \in S$, $\nu \in E$, $n\in \omega \}\cup \{z_{\delta ,j}:\delta
\in S$, $j\in \sigma \}\cup \{x_{\beta }:\beta \in \lambda \}$ and $K$ is
the subgroup with basis $\{w_{\delta ,\nu ,n}:\delta \in S$, $\nu \in E$, $%
n\in \omega \}$ where 
\begin{equation}
w_{\delta ,\nu ,n}=2y_{\delta ,\nu ,n+1}-y_{\delta ,\nu ,n}-z_{\delta ,\zeta
_{\nu }(n)}+x_{\eta _{\delta }(\nu +n)}\text{.}
\end{equation}
(If $\sigma =\aleph _{0}$, let $E=\{0\}$ and  omit $\bar{\zeta}$ and the $%
z_{\delta ,j}$.) For future reference, let $F_{\alpha }$ be the subgroup of $%
F$ generated by $\{y_{\delta ,\nu ,n}:\delta \in S\cap \alpha $, $\nu \in E$%
, $n\in \omega \}\cup \{z_{\delta ,j}:\delta \in S\cap \alpha $, $j<\sigma
\}\cup \{x_{\beta }:\beta <\alpha \}$ and for $\alpha \in S$ and $\tau
<\sigma $ let $F_{\alpha ,\tau }$ be the subgroup generated by $\{z_{\alpha
,j}:j<\tau \}$.
\end{definition}

When the cardinality of $A$ is singular, we will use a special case of a
recent result of the second author \cite{Sh667}. For convenience, we give
the statement and proof of this ``very weak diamond'' result here.

\begin{lemma}
\label{vwdmd}Assume GCH. Let $\mu $ be a singular cardinal and let $\sigma =%
\limfunc{cof}(\mu )$ and $\lambda =\mu ^{+}$. Suppose that $S$ is a
stationary subset of $\lambda $ consisting of ordinals of cofinality $\sigma 
$ and $\{\eta _{\delta }:\delta \in S\}$ is a ladder system on $S$. Then for
each $\delta \in S$ there is a sequence of sets $D^{\delta }=\left\langle
D_{\nu }^{\delta }:\nu <\sigma \right\rangle $ such that

(a) for all $\delta \in S$ and $\nu \in \sigma $, $D_{\nu }^{\delta
}\subseteq \lambda $, $\sup (D_{\nu }^{\delta })<\delta $ and $|D_{\nu
}^{\delta }|<\mu $; and

(b) for every function $h:\lambda \rightarrow \lambda $, $\{\delta \in
S:h(\eta _{\delta }(\nu ))\in D_{\nu }^{\delta }$ for all $\nu \in \sigma \}$
is stationary in $\lambda $.
\end{lemma}

\noindent \textsc{Proof}. Fix $\delta \in S$. Let $\left\langle b_{\nu
}^{\delta }:\nu <\sigma \right\rangle $ be an increasing continuous union of
subsets of $\delta $ whose union is $\delta $ and such that $\sup (b_{\nu
}^{\delta })<\delta $ and card$(b_{\nu }^{\delta })$ $<\mu $. Let $\theta
=2^{\sigma }=\sigma ^{+}$($<\mu $) and let $\left\langle g_{i}:i<\theta
\right\rangle $ be a list of all functions from $\sigma $ to $\sigma $. Also
let $\left\langle f_{\gamma }:\gamma <\lambda \right\rangle $ list all
functions from $\theta $ to $\lambda $ ($=2^{\mu }=\lambda ^{\theta }$). For
each $i\in \theta $ and $\nu \in \theta $, define $D_{\nu }^{i,\delta
}=\{f_{\gamma }(i):\gamma \in b_{g_{i}(\nu )}^{\delta }\}$.

We claim that for some $i\in \theta $, the sets $\{D^{i,\delta }=$ $%
\left\langle D_{\nu }^{i,\delta }:\nu <\sigma \right\rangle :\delta \in S\}$
will work in (b). Assuming the contrary, for each $i\in \theta $, let $%
h_{i}:\lambda \rightarrow \lambda $ be a counterexample, i.e., there is a
club $C_{i}$ in $\lambda $ such that for each $\delta \in C_{i}\cap S$,
there is $\nu \in \sigma $ such that $h_{i}(\eta _{\delta }(\nu ))\notin
D_{\nu }^{i,\delta }$.

For each $\alpha \in \lambda $, there is $h(\alpha )\in \lambda $ such that
for all $i\in \theta $, $h_{i}(\alpha )=f_{h(\alpha )}(i)$. There exists $%
\delta _{*}\in \bigcap_{i\in \theta }C_{i}\cap S$ such that for all $\alpha
<\delta _{*}$, $h(\alpha )\in \delta _{*}$. Denote $h(\eta _{\delta
_{*}}(\nu ))$ by $\gamma _{\nu }$. There exists $i_{*}\in \theta $ such that
for all $\nu <\sigma $, 
\[
g_{i_{*}}(\nu )=\min \{j<\sigma :\gamma _{\nu }\in b_{j}^{\delta _{*}^{{}}}\}%
\text{.}
\]
(Note that the right-hand side exists since $\delta _{*}=$ $\cup _{j<\sigma
}b_{j}^{\delta _{*}}$ and $\gamma _{\nu }\in \delta _{*}$.) Thus 
\[
\gamma _{\nu }\in b_{g_{i_{*}}(\nu )}^{\delta _{*}}\text{.}
\]
But then, (letting $\alpha =\eta _{\delta _{*}}(\nu )$ in the definition of $%
h$), 
\[
h_{i_{*}}(\eta _{\delta _{*}}(\nu ))=f_{h(\eta _{\delta _{*}}(\nu
))}(i_{*})=f_{\gamma _{\nu }}(i_{*})\in D_{\nu }^{i_{*},\delta _{*}}\text{.}
\]
Since this holds for all $\nu \in \sigma $, the fact that $h_{i_{*}}$ is a
counterexample implies that $\delta _{*}\notin C_{i_{*}}\cap S$. But this
contradicts the choice of $\delta _{*}$.\qed  

\medskip

\begin{theorem}
\label{nonvan1}Assume GCH. Let $\mu $ be a singular cardinal of cofinality $%
\sigma $. If $H$ is a group of cardinality $\lambda =\mu ^{+}$ built on $%
\bar{\eta}$ and $\bar{\zeta}$ as in Definition \ref{built} and $A$ is a
Whitehead group of cardinality $\mu $, then $\limfunc{Ext}(H,A)\neq 0$.
\end{theorem}

\noindent \textsc{Proof}. Let the sets $\{D^{\delta }=\left\langle D_{\nu
}^{\delta }:\nu \in \sigma \right\rangle :\delta \in S\}$ be as in Lemma \ref
{vwdmd} for this ladder system. Write $A=\bigcup_{\nu <\mu }A_{\nu }$ as in
Lemma \ref{nice2}. Without loss of generality we can assume that the
universe of $A\ $is $\mu $. 

We claim that for all $\beta <\mu $, the $2$-adic completion of $A/A_{\beta }
$ has rank $\geq \mu $ over $A/A_{\beta }$. For notational convenience we
will prove the case $\beta =0$, but the argument is the same in general
using the decomposition $A/A_{\beta }=\bigcup_{\beta \leq \alpha <\mu
}A_{\alpha }/A_{\beta }$. Since $A_{\alpha +1}/A_{\alpha }$ is $\aleph _{1}$%
-free and non-zero, there are $s_{n}^{\alpha }\in A_{\alpha +1}$ such that
the element $\Sigma _{n\in \omega }2^{n}(s_{n}^{\alpha }+A_{\alpha })$ of
the $2$-adic completion of $A_{\alpha +1}/A_{\alpha }$ is not in $A_{\alpha
+1}/A_{\alpha }$. We claim that the elements $\{\Sigma _{n\in \omega
}2^{n}s_{n}^{\alpha }:\alpha \in \mu \}$ of the $2$-adic completion of $A$
are linearly independent over $A$. Suppose not, and let 
\[
\Sigma _{i=1}^{m}k_{i}(\Sigma _{n\in \omega }2^{n}s_{n}^{\alpha (i)})=a
\]
be a counterexample; so $a\in A$; $k_{i}\in \Bbb{Z}-\{0\}$; and $\alpha
(1)<\alpha (2)<...<\alpha (m)<\mu $. Let $\gamma =\alpha (m)$ and $%
k=k_{\gamma }$. We claim that the element $k\Sigma _{n\in \omega
}2^{n}(s_{n}^{\gamma }+A_{\gamma })$ of the $2$-adic completion of $%
A_{\gamma +1}/A_{\gamma }$ belongs to $A_{\gamma +1}/A_{\gamma }$ which is a
contradiction of the choice of the $s_{n}^{\gamma }$. Since $A/A_{\gamma +1}$
is $\aleph _{1}$-free, we can write $\left\langle A_{\gamma
+1},a\right\rangle _{*}=A_{\gamma +1}\oplus C$ for some $C$, and let $%
a^{\prime }$ be the projection of $a$ on the first factor. For every $r\in
\omega $, $2^{r+1}$divides $a-$ $\Sigma _{i=1}^{m}k_{i}(\Sigma
_{n=0}^{r}2^{n}s_{n}^{\alpha (i)})$ in $A$ and hence $2^{r+1}$divides $%
a^{\prime }-$ $\Sigma _{i=1}^{m}k_{i}(\Sigma _{n=0}^{r}2^{n}s_{n}^{\alpha
(i)})$ in $A_{\gamma +1}$. But then $2^{r+1}$ divides $(a^{\prime
}+A_{\gamma })-$ $k\Sigma _{n=0}^{r}2^{n}(s_{n}^{\gamma }+A_{\gamma })$ in $%
A_{\gamma +1}/A_{\gamma }$.

Choose a strictly increasing continuous function $\xi :\sigma \rightarrow
\mu $ whose range is cofinal in $\mu $. For each $\delta \in S$ and $\nu \in
E$, there is an element $a_{\delta ,\nu }=\Sigma _{n\in \omega
}2^{n}(a(\delta ,\nu ,n)+A_{\xi (\nu )+1})$ in the $2$-adic completion of $%
A/A_{\xi (\nu )+1}$ which is not in the subgroup generated by $A/A_{\xi (\nu
)+1}$ and the $2$-adic completion of $\{d+A_{\xi (\nu )+1}:d\in D_{\nu
}^{\delta }\cap A\}$. (Note that the latter has cardinality $<$ $\mu $ since 
$|D_{\nu }^{\delta }|^{\aleph _{0}}<\mu $ by the GCH.)

Now define $\psi :K\rightarrow A$ such that $\psi (w_{\delta ,\nu
,n})=a(\delta ,\nu ,n)$. We claim that $\psi $ does not extend to a
homomorphism $\varphi :F\rightarrow A$. Suppose, to the contrary, that it
does. Then by Lemma \ref{vwdmd}, there is $\delta \in S$ such that $\varphi
(x_{\eta _{\delta }(\nu )})\in D_{\nu }^{\delta }$ for all $\nu \in \sigma $%
. Now there exists $\nu \in E$ such that $\varphi (z_{\delta ,j})\in A_{\xi
(\nu )}$ for all $j<\nu $. We will contradict the choice of $a_{\delta ,\nu
} $ for this $\delta $ and $\nu $.

We work in $A/A_{\xi (\nu )+1}$. Let $c_{n}=\varphi (y_{\delta ,\nu
,n})+A_{\xi (\nu )+1}$, $d_{n}=\varphi (x_{\eta _{\delta }(\nu +n)})+A_{\xi
(\nu )+1}$. Then by applying $\varphi $ to the equations (\ref{eqnsG}) and
since $\varphi (z_{\delta ,j})\in A_{\xi (\nu )}$ for all $j<\nu $ we have
that for all $n\in \omega $, 
\[
a(\delta ,\nu ,n)+A_{\xi (\nu )+1}=2c_{n+1}-c_{n}+d_{n}\text{.}
\]
It follows that $a_{\delta ,\nu }=c_{0}+\Sigma _{n\in \omega }2^{n}d_{n}$ is
in the subgroup generated by $A/A_{\xi (\nu )+1}$ and the $2$-adic
completion of $\{d+A_{\xi (\nu )+1}:d\in D_{\nu }^{\delta }\cap A\}$, which
contradicts the choice of $a_{\delta ,\nu }$. \qed  

\medskip

We now turn to the cases when the cardinality of $A$ is a successor
cardinal. Though the two arguments could be combined into one, following the
argument in Theorem \ref{nonvan3}, we prefer to introduce the method with
the somewhat simpler argument for the successor of regular case. The
following lemma is easy to confirm:

\begin{lemma}
\label{2adic}Suppose that $L^{\prime }$ is a free subgroup of $L$ such that $%
L/L^{\prime }$ is $\aleph _{1}$-free. If $\{t_{n}:n\in \omega \}$ is a basis
of a summand of $L^{\prime }$, then $\sum_{n\in \omega }2^{n}t_{n}$ is an
element of the $2$-adic completion of $L$ which does not belong to $L$. In
other words, the system of equations 
\[
2y_{n+1}=y_{n}-t_{n}
\]
in the unknowns $y_{n}$ ($n\in \omega $) does not have a solution in $L$. 
\qed  
\end{lemma}

\medskip

\begin{theorem}
\label{nonvan2}Assume GCH. Let $\lambda =\mu ^{+}$ where $\mu $ is a regular
cardinal. Suppose $H$ is built on $\bar{\eta}=\{\eta _{\delta }:\delta \in
S\}$ and $\bar{\zeta}$ $=\{\zeta _{\nu }:\nu \in E\}$ as in Definition \ref
{built}. Suppose also, for $\mu >\aleph _{0}$, that $\diamondsuit _{\mu
}(E^{\prime })$ holds for all stationary subsets $E^{\prime }$ of $E$. If $A$
is a Whitehead group of cardinality $\lambda =\mu ^{+}$, then $\limfunc{Ext}%
(H,A)\neq 0$.
\end{theorem}

\noindent \textsc{Proof}. Let $A=\bigcup_{\alpha <\lambda }A_{\alpha }$ and $%
S(A)$ be as in Lemma \ref{nice}. Note that we make no assumption about the
relation of $S$ and $S(A)$; maybe $S\cap S(A)=\emptyset $. Without loss of
generality, for all $\delta \in S(A)$, $A_{\delta +1}/A_{\delta }$ is
Whitehead of rank $\mu $ and $A/A_{\delta +1}$ is Whitehead. Assume $\mu
>\aleph _{0}$; the proof for $\aleph _{0}$ is simpler. For each $\delta
<\lambda $, write $A_{\delta }$ as the union of a continuous chain of
subgroups of cardinality $<$ $\mu $: $A_{\delta }=\bigcup_{\nu <\mu
}B_{\delta ,\nu }$. For $\delta \in S(A)$, since $\diamondsuit _{\mu }(E)$
holds, we can assume that the set of $\nu \in E$ such that $A_{\delta
+1}/(A_{\delta }+B_{\delta +1,\nu })$ is locally Whitehead is stationary;
for such $\nu $, the quotient is then strongly $\aleph _{1}$-free since CH
holds. Thus for $\nu $ in a stationary subset $E_{\delta }$ of $E$ we can
assume that $A_{\delta }+B_{\delta +1,\nu +1}/A_{\delta }+B_{\delta +1,\nu }$
is free of rank $\aleph _{0}$ and $A_{\delta +1}/A_{\delta }+B_{\delta
+1,\nu +1}$ is $\aleph _{1}$-free. Say $\{t_{\delta ,\nu ,n}+A_{\delta
}+B_{\delta +1,\nu }:n\in \omega \}$ is a basis of $A_{\delta }+B_{\delta
+1,\nu +1}/A_{\delta }+B_{\delta +1,\nu }$.

For each $\delta _{1}\in S$, let $\delta _{1}^{+}$ be the least member of $%
S(A)$ which is $\geq \delta _{1}$. Define 
\[
\psi (w_{\delta _{1},\nu ,n})=t_{\delta _{1}^{+},\nu ,n}
\]
for all $w_{\delta _{1},\nu ,n}\in K$ if $\nu \in E_{\delta _{1}^{+}}$. We
claim that $\psi $ does not extend to $\varphi :F\rightarrow A$. Suppose to
the contrary that it does. Let $M=\varphi [F]$, $M_{\alpha }=\varphi
[F_{\alpha }]$, $M_{\alpha ,\tau }=\varphi [F_{\alpha ,\tau }]$. Then there
is a club $C$ in $\lambda $ such that for $\alpha \in C$, $M_{\alpha
}\subseteq A_{\alpha }$. Fix $\delta _{1}$ in $C\cap S$. Let $\delta $ be $%
\delta _{1}^{+}$ and choose $\gamma \in C$ such that $\gamma >\delta $.
There is a club $C^{\prime }$ in $\mu $ such that for $\nu \in C^{\prime }$, 
$M_{\delta _{1},\nu }\subseteq B_{\gamma ,\nu }$ and $A_{\delta +1}\cap
B_{\gamma ,\nu }\subseteq B_{\delta +1,\nu }$. Since $\diamondsuit _{\mu
}(E_{\delta })$ holds, there is, by Lemma\ref{nice},  $\nu \in E_{\delta
}\cap C^{\prime }$ such that $A_{\gamma }/(A_{\delta +1}+B_{\gamma ,\nu })$
is locally Whitehead, and hence $\aleph _{1}$-free. We will obtain a
contradiction of Lemma \ref{2adic} with $L=A_{\gamma }/A_{\delta }+B_{\gamma
,\nu }$ and $L^{\prime }=(B_{\delta +1,\nu +1}+A_{\delta }+B_{\gamma ,\nu
})/A_{\delta }+B_{\gamma ,\nu }$ and $t_{n}=t_{\delta ,\nu ,n}+A_{\delta
}+B_{\gamma ,\nu }$. Notice that modulo $A_{\delta }+B_{\gamma ,\nu }$ we
have 
\[
2\varphi (y_{\delta _{1},\nu ,n+1})=\varphi (y_{\delta _{1},\nu
,n})-t_{\delta ,\nu ,n}
\]
for all $n\in \omega $ since $\varphi (x_{\eta _{\delta _{1}}(\nu +n)})\in
A_{\delta }$ and $\varphi (z_{\delta _{1},\zeta _{\nu }(n)})\in B_{\gamma
,\nu }$. Moreover, $\{t_{n}:n\in \omega \}$ is a basis of a summand of $%
L^{\prime }$ since $L^{\prime }$ is naturally isomorphic to $A_{\delta }$ $%
+B_{\delta +1,\nu +1}/A_{\delta }+(B_{\gamma ,\nu }\cap (A_{\delta
}+B_{\delta +1,\nu +1}))$ and the latter has a natural epimorphism onto $%
A_{\delta }+B_{\delta +1,\nu +1}/A_{\delta }+B_{\delta +1,\nu }$ which is
free on the basis $\{t_{\delta ,\nu ,n}+A_{\delta }+B_{\delta +1,\nu }:n\in
\omega \}$. It remains to show that $L/L^{\prime }$ is $\aleph _{1}$-free.
Now 
\[
0\rightarrow (A_{\delta +1}+B_{\gamma ,\nu }))(B_{\delta +1,\nu
+1}+A_{\delta }+B_{\gamma ,\nu })\rightarrow L/L^{\prime }\rightarrow
A_{\gamma }/(A_{\delta +1}+B_{\gamma ,\nu })\rightarrow 0
\]
is exact and $A_{\gamma }/(A_{\delta +1}+B_{\gamma ,\nu })$ is $\aleph _{1}$%
-free by choice of $\nu $, so it suffices to show that $(A_{\delta
+1}+B_{\gamma ,\nu })/(B_{\delta +1,\nu +1}+A_{\delta }+B_{\gamma ,\nu })$
is $\aleph _{1}$-free. But this is isomorphic to $A_{\delta +1}/((A_{\delta
}+B_{\delta +1,\nu +1})+(A_{\delta +1}\cap B_{\gamma ,\nu }))$, which (since 
$A_{\delta +1}\cap B_{\gamma ,\nu }\subseteq B_{\delta +1,\nu }\subseteq
B_{\delta +1,\nu +1}$) equals $A_{\delta +1}/(A_{\delta }+B_{\delta +1,\nu
+1})$, which was chosen $\aleph _{1}$-free. \qed  

\begin{theorem}
\label{nonvan3}Assume GCH. Let $\lambda =\mu ^{+}$ where $\mu $ is a
singular cardinal of cofinality $\sigma <\mu $. Suppose $H$ is built on $%
\bar{\eta}=\{\eta _{\delta }:\delta \in S\}$ and $\bar{\zeta}$ $=\{\zeta
_{\nu }:\nu \in E\}$ as in Definition \ref{built}. Suppose also that $%
\diamondsuit _{\lambda }(Y)$ holds for some subset $Y$ of $\lambda $
consisting of limit ordinals of cofinality $\sigma $ and that, if $\sigma
>\aleph _{0}$, $\diamondsuit _{\sigma }(E)$ holds. If $A$ is a Whitehead
group of cardinality $\lambda =\mu ^{+}$, then $\limfunc{Ext}(H,A)\neq 0$.
\end{theorem}

\noindent \textsc{Proof}. Without loss of generality, for all $\delta \in
S(A)$, $A_{\delta +1}/A_{\delta }$ is Whitehead of rank $\mu $. For each $%
\delta \in S$, choose a strictly increasing continuous sequence $%
\left\langle \xi _{\delta ,\nu }:\nu \leq \sigma \right\rangle $ of elements
of $S(A)$ such that $\xi _{\delta ,0}\geq \delta +1$. (This is possible
because, by Lemma \ref{nice}, $Y\cap S(A)$ is stationary so we can choose $%
\xi _{\delta ,\sigma }$ to be an element of $Y\cap S(A)\cap \overline{%
(S(A)\cap (\delta ,\lambda ))}$ where $\overline{(S(A)\cap (\delta ,\lambda
))}$ is the closure of $\{\alpha \in S(A):\alpha >\delta \}$.) Let $%
B_{\delta +1,\nu }=A_{\xi _{\delta ,\nu }}$. (Note the difference from the
last proof.) We can then modify the sequence so that $B_{\delta +1,\nu
+1}/B_{\delta +1,\nu }$ is free on a countable set $\{t_{\delta ,\nu
,n}+B_{\delta +1,\nu }\}$ when $\nu \in E$.

For each $\delta _{1}\in S$, let $\delta _{1}^{+}$ be the least member of $%
S(A)$ which is $\geq \delta _{1}$. Define 
\[
\psi (w_{\delta _{1},\nu ,n})=t_{\delta _{1}^{+},\nu ,n} 
\]
for all $w_{\delta _{1},\nu ,n}\in K$. We claim that $\psi $ does not extend
to $\varphi :F\rightarrow A$. Suppose to the contrary that it does. As
before, let $M=\varphi [F]$, $M_{\alpha }=\varphi [F_{\alpha }]$, $M_{\alpha
,\tau }=\varphi [F_{\alpha ,\tau }]$ and let $C$ be a club such that for $%
\alpha \in C$, $M_{\alpha }\subseteq A_{\alpha }$. Fix $\delta _{1}$ in $%
C\cap S$. Let $\delta $ be $\delta _{1}^{+}$ and choose $\gamma \in C$ such
that $\gamma >\delta $.

Let $N=\bigcup_{\nu <\sigma }N_{\nu }$ be the continuous union of a chain of
elementary submodels of $H(\chi )$ for large enough $\chi $ such that each $%
N_{\nu }$ has cardinality $<\sigma $, $N_{\nu }\in N_{\nu +1}$ and such that 
$\delta $, $\sigma $, $A$, $\{\varphi (z_{\delta _{1},\nu }):\nu <\sigma \}$%
, $\{\varphi (x_{\eta _{\delta _{1}}(\nu +n)}):\nu <\sigma \}$ (for each $%
n\in \omega $), $\{t_{\delta ,\nu ,n}:\nu <\sigma ,n\in \omega \}$ and $%
\{\xi _{\delta ,\nu }:\nu \leq \sigma \}$ all belong to $N_{0}$ and 
\[
\{\varphi (z_{\delta _{1},j}):j<\sigma \}\cup \{\varphi (x_{\eta _{\delta
_{1}}(\nu )}):\nu <\sigma \}\cup \{t_{\delta ,\nu ,n}:\nu <\sigma ,n\in
\omega \}\cup \sigma \subseteq N\text{. }
\]
Moreover, by intersecting with a club, we can assume that for all $\nu $, $%
N_{\nu }\cap \sigma =\nu $ and $N_{\nu }\cap B_{\delta +1,\sigma }\subseteq
B_{\delta +1,\nu }$ and hence $\{\xi _{\delta ,j}:j<\nu \}$, $\{\varphi
(z_{\delta _{1},j}):j<\nu \}$, $\{t_{\delta ,j,n}:j<\nu ,n\in \omega \}$,
and $\{\varphi (x_{\eta _{\delta _{1}}(j+n)}):j<\nu \}$ (for all $n\in
\omega $) are all subsets of $N_{\nu }$. We claim that there is a $\nu \in E$
such that $A/(B_{\delta +1,\sigma }+(N_{\nu }\cap A))$ is Whitehead, and
hence $\aleph _{1}$-free. Assuming this for the moment, we show how to
obtain a contradiction of Lemma \ref{2adic} with $L=(N\cap A)/(N\cap
A_{\delta })+(N_{\nu }\cap A)$, $L^{\prime }=((N\cap B_{\delta +1,\nu
+1})+(N_{\nu }\cap A))/(N\cap A_{\delta })+(N_{\nu }\cap A)$ and $%
t_{n}=t_{\delta ,\nu ,n}+(N_{\nu }\cap A)$. Notice that for all $n\in \omega 
$, $\varphi (x_{\eta _{\delta _{1}}(\nu +n)})\in (N\cap A_{\delta })$ and $%
\varphi (z_{\delta _{1},\zeta _{\nu }(n)})\in N_{\nu }$. Moreover, $%
\{t_{n}:n\in \omega \}$ is a basis of a summand of $L^{\prime }$ because $%
L^{\prime }$ is isomorphic to $(N\cap B_{\delta +1,\nu +1})/(N\cap A_{\delta
})+(N_{\nu }\cap B_{\delta +1,\nu })$ and the latter has epimorphic image $%
(N\cap B_{\delta +1,\nu +1})/(N\cap B_{\delta +1,\nu })$ which is free on
the basis $\{t_{\delta ,\nu ,n}+(N\cap B_{\delta +1,\nu }):n\in \omega \}$.
To see that $L/L^{\prime }$ is $\aleph _{1}$-free, use the short exact
sequence

\begin{quotation}
$0 \rightarrow ((N\cap B_{\delta +1,\sigma })+(N_{\nu }\cap A))/((N\cap
B_{\delta +1,\nu +1})+(N_{\nu }\cap A))\rightarrow
\newline 
L/L^{\prime }\rightarrow $
$(N\cap A)/((N\cap B_{\delta +1,\sigma })+(N_{\nu }\cap A)) \rightarrow 0%
\text{.} $
\end{quotation}

The last term is $\aleph _{1}$-free by choice of $\nu $ and since $N$ is an
elementary submodel of $H(\chi )$. Moreover, $((N\cap B_{\delta +1,\sigma
})+(N_{\nu }\cap A))/((N\cap B_{\delta +1,\nu +1})+(N_{\nu }\cap A))$ is
isomorphic to $(N\cap B_{\delta +1,\sigma })/(N\cap B_{\delta +1,\nu +1})$
(since $N_{\nu }\cap B_{\delta +1,\sigma }\subseteq B_{\delta +1,\nu }$) and
thus is $\aleph _{1}$-free since $A/B_{\delta +1,\nu +1}$ is $\aleph _{1}$%
-free.

It remains to show that there is a $\nu \in E$ such that $A/(B_{\delta
+1,\sigma }+(N_{\nu }\cap A))$ is Whitehead. If not, then for all $\nu \in E$%
, $(B_{\delta +1,\sigma }+(N_{\nu +1}\cap A))/(B_{\delta +1,\sigma }+(N_{\nu
}\cap A))$ is not Whitehead, since $B_{\delta +1,\sigma }$, $A$ and $N_{\nu
} $ belong to the elementary submodels $N_{\nu +1}$ and $N$. But then $%
\diamondsuit _{\sigma }(E)$ implies that $\bigcup_{\nu <\sigma }(B_{\delta
+1,\sigma }+(N_{\nu }\cap A))/B_{\delta +1,\sigma }$ is a group of
cardinality $\sigma $ which is not a Whitehead group, contradicting the fact
that $A/B_{\delta +1,\sigma }=A/A_{\xi _{\delta ,\sigma }}$ is locally
Whitehead. \qed  

\section{Whitehead groups by uniformization}

We present a special case of a theorem of Shelah and Str\"{u}ngmann \cite{SS}%
.

\begin{theorem}
\label{unif}Suppose that $H$ is built from $\bar{\eta}$ and $\bar{\zeta}$ as
in Definition \ref{built} and that $E$ is a non-reflecting subset of $\sigma 
$. Then $H$ is a Whitehead group if $\bar{\eta}$ satisfies $\omega $%
-uniformization, that is, for every family of functions $\{c_{\delta
}:\sigma \rightarrow \omega :\delta \in S\}$, there is a pair $(f,f^{*})$
where $f:\lambda \rightarrow \omega $ and $f^{*}:S\rightarrow \sigma $ such
that $f(\eta _{\delta }(\nu ))=c_{\delta }(\nu )$ whenever $f^{*}(\delta
)\leq \nu <\sigma $.
\end{theorem}

\noindent \textsc{Proof}. We assume $\sigma >\aleph _{0}$ since this is
known otherwise (cf. \cite{Sh80}, \cite{T}). If $F$ and $K$ are as in
Definition \ref{built}, it suffices to show that every homomorphism $\psi
:K\rightarrow \Bbb{Z}$ extends to a homomorphism $\varphi :F\rightarrow \Bbb{%
Z}$. Given $\psi $, define $c_{\delta }(\nu +n)=\psi (w_{\delta ,\nu ,n})$
for $\nu \in E$, and arbitrary otherwise. Let $(f,f^{*})$ be the
uniformizing pair. Define $\varphi (x_{\beta })=f(\beta )$. For each $\delta
\in S$ we must still define $\varphi (y_{\delta ,\nu ,n})$ and $\varphi
(z_{\delta ,j})$ for $\nu ,j\in \sigma $ and $n\in \omega $. Fix $\delta $
and let $\rho =f^{*}(\delta )$; without loss of generality $\rho \notin E$.
Let $F^{\prime }$ (resp. $F_{\rho }^{\prime }$) be the subgroup of $F$
generated by $\{y_{\delta ,\nu ,n}:$ $\nu \in E$, $n\in \omega \}\cup
\{z_{\delta ,j}:$ $j<\sigma \}\cup \{x_{\beta }:\beta <\delta \}$ (resp. by $%
\{y_{\delta ,\nu ,n}:$ $\nu \in E\cap \rho $, $n\in \omega \}\cup
\{z_{\delta ,j}:$ $j<\rho \}\cup \{x_{\beta }:\beta <\delta \}$ ) and $%
K^{\prime }$ (resp., $K_{\rho }^{\prime }$) the subgroup generated by $%
\{w_{\delta ,\nu ,n}:$ $\nu \in E$, $n\in \omega \}\cup \{x_{\beta }:\beta
<\delta \}$ (resp., by $\{w_{\delta ,\nu ,n}:$ $\nu \in E\cap \rho $, $n\in
\omega \}\cup \{x_{\beta }:\beta <\rho \}$). Then $F^{\prime }/K^{\prime }$
is $\sigma $-free since $E$ is non-reflecting, so $K_{\rho }^{\prime }$ is a
summand of $F_{\rho }^{\prime }$; then it is easy to extend $\psi
\upharpoonright \{w_{\delta ,\nu ,n}:$ $\nu \in E\cap \rho $, $n\in \omega
\}+\varphi \upharpoonright \{x_{\beta }:\beta <\rho \}$ to $\varphi :F_{\rho
}^{\prime }\rightarrow \Bbb{Z}$. For $\nu \in E$ with $\nu >\rho $ we have $%
\varphi (x_{\eta _{\delta }(\nu +n)})=\psi (w_{\delta ,\nu ,n})$ for all $%
n\in \omega $. For some $m_{\nu }$, $\zeta _{\nu }(n)\geq \rho $ when $n\geq
m_{\nu }$. Then we can satisfy the equations 
\[
\psi (w_{\delta ,\nu ,n})=2\varphi (y_{\delta ,\nu ,n+1})-\varphi (y_{\delta
,\nu ,n})-\varphi (z_{\delta ,\zeta _{\nu }(n)})+\varphi (x_{\eta _{\delta
}(\nu +n)})
\]
by setting $\varphi (y_{\delta ,\nu ,n})=0=\varphi (z_{\delta ,\zeta _{\nu
}(n)})$ for $n\geq m_{\nu }$ and defining $\varphi (y_{\delta ,\nu ,n})$ by
downward induction for $n<m_{\nu }$. \qed  

\medskip

Finally we can put the pieces together to prove:

\begin{theorem}
\label{setcons}There is a model of ZFC + GCH such that for every Whitehead
group $A$ of infinite rank, there is a Whitehead group $H_{A}$ of
cardinality $\leq |A|^{+}$ such that $\limfunc{Ext}(H_{A},A)\neq 0$.
\end{theorem}

\noindent \textsc{Proof}. By standard forcing methods (cf. \cite{ES}) there
is a model of ZFC + GCH such that:

\begin{quotation}
 \ 

(i) for every infinite successor cardinal $\lambda =\mu ^{+}$ there is a
stationary subset $S$ of $S_{\limfunc{cf}(\mu )}^{\lambda }$ with a ladder
system $\bar{\eta}=\{\eta _{\delta }:\delta \in S\}$ which satisfies $\omega 
$-uniformization (or even $\kappa $-uniformization for $\kappa <\mu $);

(ii) for every infinite successor cardinal $\lambda =\mu ^{+}$ there is a
stationary subset $Y$ of $S_{\limfunc{cf}(\mu )}^{\lambda }$ such that $%
\diamondsuit _{\lambda }(Y)$ holds;

(iii) for every regular uncountable cardinal $\sigma $, there is a
non-reflecting stationary subset $E$ of $S_{\omega }^{\sigma }$ such that $%
\diamondsuit _{\sigma }(E^{\prime })$ holds for every stationary subset $%
E^{\prime }$ of $E$;

(iv) there is a tree-like ladder system on a stationary subset of $\omega
_{1}$ which satisfies $2$-uniformization but not $\omega $-uniformization.
\end{quotation}

(In fact, we can get more: we can strengthen (i) and (ii) to the following:
for every infinite successor cardinal $\lambda =\mu ^{+}$ there is a normal
ideal $I_{\lambda }$ containing the non-stationary ideal such that for every 
$S\in I_{\lambda }$, $S-S_{\limfunc{cf}(\mu )}^{\lambda }$ is
non-stationary, and there exists a stationary $S^{\prime }\in I_{\lambda }$
disjoint from $S$; moreover, for every $S\in I_{\lambda }$, there is a
ladder system $\bar{\eta}=\{\eta _{\delta }:\delta \in S\}$ which satisfies $%
\omega $-uniformization and for every $S\notin I_{\lambda }$, $\diamondsuit
_{\lambda }(S)$ holds.)

We work in this model. Let $A$ be a Whitehead group of infinite rank. If the
rank of $A$ is $\aleph _{0}$, then $A$ is isomorphic to $\Bbb{Z}^{(\omega )}$
and it is well-known (cf. \cite{Sh80}, \cite[XII.3.]{EM}) that (iv) implies
that there is a Whitehead group $H$ which is not $\aleph _{1}$-coseparable,
i.e., $\limfunc{Ext}(H,\Bbb{Z}^{(\omega )})\neq 0$. If the cardinality of $A$
is either singular or a successor cardinal, then for $\lambda =|A|$ if $|A|$
is regular, or $\lambda =|A|^{+}$ if $|A|$ is singular, the properties (i),
(ii) and (iii) allow us to build a group $H_{A} $ of cardinality $\lambda $
as in Definition \ref{built}, which is Whitehead by Theorem \ref{unif} and
such that by Theorem \ref{nonvan1}, \ref{nonvan2} or \ref{nonvan3}, $%
\limfunc{Ext}(H_{A,}A)\neq 0$.

It is also consistent to assume that there are no regular limit (i.e.
inaccessible) cardinals, in which case we have covered all possibilities for
the cardinality of $A$ and we are done. Another approach is to allow
inaccessible cardinals but force the model to satisfy in addition:

\begin{quote}
(v) for every inaccessible cardinal $\lambda $ there is a stationary subset $%
S$ of $S_{\aleph _{0}}^{\lambda }$ with a ladder system $\bar{\eta}=\{\eta
_{\delta }:\delta \in S\}$ which satisfies $\omega $-uniformization;
moreover $\diamondsuit _{\lambda }$ holds.
\end{quote}

\noindent As in Lemma \ref{nice}, one can show that $S(A)$ is stationary and
then the proof is similar to that in Theorem \ref{nonvan2}. \qed

\end{document}